\input amstex

\documentstyle{amsppt}
\document
\topmatter
\title
On $\Cal A$-tensors in  Riemannian geometry.
\endtitle
\author
Wlodzimierz Jelonek
\endauthor
\abstract{We present examples, both compact and non-compact complete,
 of locally non-homogeneous proper $\Cal A$-manifolds.}
 \endabstract

\endtopmatter
\define\G{\Gamma}

\define\bt{\bar\theta}
\define\n{\nabla}
\define\bx{\bar\xi}
\define\w{\wedge}
\define\k{\diamondsuit}
\define\th{\theta}
\define\p{\parallel}
\define\a{\alpha}
\define\g{\gamma}
\define\lb{\lambda}
\define\ra{\rangle}
\define\la{\langle}
\define\A{\Cal A}
\define\1{D_{\lb}}
\define\2{D_{\mu}}
\define\bO{\bar{\Omega}}
\define\wO{\tilde{\Omega}}
\define\wT{\tilde{T}}
\define\0{\Omega}
\define\wn{\tilde\nabla}
\par

\medskip
{\bf 0. Introduction.} A.Gray in the paper [7] introduced
the notion of $\Cal A$ - manifolds. An $\Cal A$-manifold it is a
Riemannian manifold $(M,g)$ whose Ricci tensor $\rho$ satisfies the
following condition : $\n \rho (X,X,X)=0$ for all $X\in TM$
where $\n$ is the Levi-Civita connection of the metric $g$.
J.E. D'Atri and H.K. Nickerson (see [3],[4], [5])
proved that every naturally reductive
homogeneous Riemannian space is a manifold with volume preserving
local geodesic symmetries.
It is also known that every Riemannian manifold with
volume preserving geodesic symmetries  belongs to
the class of $\A$-manifolds (see [14]).
Also in [7] and [12] there are many examples of $\A$-manifolds.
All of them are (locally) homogeneous.
The aim of our paper is to investigate symmetric tensors $S\in$ End$(TM)$
on a Riemannian manifold $(M, g)$ satisfying an additional
condition: $\n\Phi (X,X,X)=0$ for all $X\in TM$, where
a tensor $\Phi$ is defined as follows: $\Phi(X,Y)=g(SX,Y)$.
We construct proper $\A$-manifolds on $S^1$-bundles
over K\"ahler-Einstein manifolds and give explicit examples of
compact locally non-homogoneous
proper $\A$-manifolds (of cohomogeneity $d$ for arbitrary $d\in \Bbb{N}$)
using N. Koiso and Y. Sakane examples of (locally) non-homogeneous
K\"ahler-Einstein
manifolds (see [9], [10]).   We give explicit
examples of locally non-homogeneous
proper complete $\A$-manifolds using E.Calabi's examples of
locally non-homogeneous
K\"ahler-Einstein complete manifolds (see [2], [16]).
In this way we give an answer to the problem
(p.451, 16.56 (i)) in the book [1].

\par
\bigskip
{\bf 1. Preleminaries.} We use the notation as in [6]. Let
$(M,g)$ be a smooth connected Riemannian manifold. Abusing the notation
we shall write $\la X,Y\ra=g(X,Y)$. For a tensor $T(X_1,X_2,..,X_k)$ we define
a tensor $\n T(X_0,X_1, ..,X_k)$ by $\n T(X_0,X_1, ..,X_k)=\n_{X_0}T(X_1,
..,X_k)$.
 By an ${\Cal A}$-tensor on $M$
we mean  an endomorphism $S\in$ End$(TM)$ satisfying the
following conditions:

(1.1)  $ \la SX,Y\ra=\la X,SY\ra$  for all $X,Y\in TM$

(1.2)  $\la\n S(X,X),X\ra=0$ for all $X \in TM$.

We also shall write $S\in \A$ if $S$ is an $\A$-tensor. We call $S$
 a proper $\A$-tensor if $\n S\ne 0$.
We  denote by $\Phi$ a tensor
defined by $\Phi(X,Y)=\la SX,Y\ra.$

We start with :
\par
\medskip
  {\bf Proposition 1.1.} {\it The following conditions are equivalent:}

  (a) {\it A tensor} $S$ {\it is an} $\Cal A${\it-tensor on} $(M,g)$;

  (b) {\it For every geodesic } $\gamma$ {\it on} $(M,g)$ {\it
    the function}
  $\Phi (\g '(t),\g '(t))$ {\it is constant on }  dom$\g$.

  (c)  {\it A  condition}
  $$ \n_X \Phi(Y,Z)+ \n_Z \Phi(X,Y)+\n_Y \Phi(Z,X)=0\tag A$$
  {\it is satisfied for all } $X,Y,Z\in \frak X(M)$.
\par
\medskip

  {\it Proof:} By using polarization it is easy to see that (a) is
  equivalent to (c). Let now  $X \in T_{x_0}M$ be any vector from $TM$
  and $\g$ be a geodesic satisfying an initial condition $\g '(0)=X$.
  Then
  $$ \frac d{dt}\Phi(\g '(t),\g '(t))=\n_{\g '(t)}\Phi(\g '(t),\g '(t)).
  \tag 1.3$$
  Hence $ \frac d{dt}\Phi(\g '(t),\g '(t))_{t=0}=\n\Phi(X,X,X).$
  The equivalence (a) $\Leftrightarrow$ (b) follows immediately
  from the above relations.$\k$
  \par
  \medskip
   Define as in [6] the integer-valued function  $E_S(x)=($ the number
   of distinct eigenvalues of $S_x$) and set $M_S=\{ x\in M:E_S$ is
   constant in a neighbourhood of $x\}$. The set $M_S$ is open and dense
   in $M$ and the eigenvalues $\lb_i$ of $S$ are distinct and smooth
   in each component $U$ of $M_S$. The eigenspaces $D_{\lb}=\ker (S-\lb I)$
    form smooth distributions in each component $U$ of $M_S$. By
    $\n f$ we denote the gradient of a function $f$ ($\la\n f,X\ra=df(X)$)
    and by $\G(D_{\lb})$  (resp. by $\frak X(U)$) the set of all local
    sections of the bundle
    $D_{\lb}$ (resp. of all local vector fields on $U$). Let us note that if $\lb\ne \mu$ are eigenvalues
    of $S$ then $\1$ is orthogonal to $\2$.
    \par
    \medskip

   {\bf Theorem 1.2.} {\it Let } $S$ {\it be an}
   $\Cal A$-{\it tensor on} $M$
   {\it and U be a component of } $M_S$ {\it and}
   $\lb_1,\lb_2,..,\lb_k \in C^{\infty}(U)$ {\it be eigenfunctions of }
   $S$. {\it Then for all} $X\in D_{\lb_i}$ {\it we have}
   $$ \n S(X,X)=-\frac12 \n \lb_i\parallel X \parallel ^2\tag 1.4$$
   {\it and}  $D_{\lb_i}\subset \ker d\lb_i$. { \it If } $i\ne j$
   {\it and} $X\in \G(D_{\lb_i}),Y \in \G(D_{\lb_j})$
   {\it then}

  $$ \la\n_X X ,Y\ra=\frac12\frac{Y\lb_i}{\lb_j-\lb_i}\parallel X \parallel ^2.
   \tag 1.5$$
\par
\medskip
   {\it Proof.} Let $X \in \G(D_{\lb_i})$ and $Y \in \frak X(U).$
   Then we have $SX= \lb_iX$ and
   $$  \n S(Y,X)+(S-\lb_iI)(\n_YX)=(Y\lb_i)X \tag 1.6$$
   and consequently,
   $$  \la\n S(Y,X),X\ra=(Y\lb_i)\p X\p^2. \tag 1.7$$
   Taking $Y=X$ in (1.7) we obtain by (1.2) $0=X\lb_i\p X\p^2$.
   Hence $D_{\lb_i}\subset \ker d\lb_i$. Thus from (1.6) it follows that
   $ \n S(X,X)=(\lb_iI-S)(\n_XX)$. Condition  (A) imply
   $ \la\n S(X,Y),Z\ra+ \la\n S(Z,X),Y\ra+\la\n S(Y,Z),X\ra=0$
   hence
    $$ 2\la\n S(X,X),Y\ra+ \la\n S(Y,X),X\ra=0.\tag 1.8 $$
    Thus, (1.8) yields
    $ Y\lb_i\p X\p^2 + 2\la\n S(X,X),Y\ra=0$. Consequently,
    $\n S(X,X)=-\frac12 \n \lb_i\p X\p^2$.
     Let now $Y\in \G(D_{\lb_j}).$
    Then we have
    $$  \n S(X,Y)+(S-\lb_jI)(\n_XY)=(X\lb_j)Y. \tag 1.9$$
    It is also clear that $\la\n S(X,X),Y\ra=\la\n S(X,Y),X\ra=
    (\lb_j-\lb_i)\la\n_XY,X\ra$. Thus,
    $$ Y\lb_i\p X \p^2=-2(\lb_j-\lb_i)\la\n_XY,X\ra
    =2(\lb_j-\lb_i)\la Y,\n_XX\ra$$
    and (1.5) holds.$\k$

    \par
    \medskip
    {\it Corrolary 1.3.} Let $S,U,\lb_1,\lb_2,..,\lb_k$ be as above and
    $i\in \{1,2,..,k\}$. Then
    the following conditions are equivalent:

    (a)  For all $X\in \G(D_{\lb_i}) \ \ \n_XX \in D_{\lb_i}.$

    (b)  For all $X,Y\in \G(D_{\lb_i}) \ \ \n_XY+\n_YX \in D_{\lb_i}.$

    (c)  For all $X\in \G(D_{\lb_i}) \ \ \n S(X,X) =0.$

    (d)  For all $X,Y\in \G(D_{\lb_i}) \ \ \n S(X,Y)+\n S(Y,X) =0.$

    (e)  $\lb_i$ is a constant eigenvalue of $S$.
    \par
    \medskip
    Let us note that if $X,Y\in \G(D_{\lb_i})$ then
    ($D_{\lb_i}\subset\ker d\lb_i $ !):
    $$\n S(X,Y)-\n S(Y,X)=(\lb_iI-S)([X,Y]) \tag 1.10$$
    hence a distribution $D_{\lb_i}$ is integrable if and only if
    $\n S(X,Y)$ $=\n S(Y,X)$ for all $X,Y\in \G(D_{\lb_i}).$  Consequently,
    we obtain
    \par
    \medskip
    {\it Corollary 1.4.} Let $\lb_i \in C^{\infty}(U)$ be an
    eigenvalue of $\Cal A$-tensor $S$. Then on $U$ the following conditions are
    equivalent:

    (a) $D_{\lb_i}$ is integrable and $\lb_i$ is constant.

    (b) For all $X,Y\in \G(D_{\lb_i})  \ \ \n S(X,Y)=0.$

    (c) $D_{\lb_i}$ is autoparallel.
    \par
    \medskip
    {\it Proof.}  It follows from (1.4),(1.10), Corollary 1.3 and the
    relation $\n_XY=\n_YX+[X,Y].$ $\k$

    \par
    \medskip

\bigskip

{\bf 2. $\Cal A$-tensors with two constant eigenvalues.} In this section
we shall characterize  $\Cal A$-tensors with two constant eigenvalues.
We start with:
\par
\medskip
{\bf Theorem 2.1.} {\it Let S be an} $\Cal A$-{\it tensor on } $(M,g)$
{\it with exactly two eigenvalues} $\lb,\mu$ {\it  and a constant trace.
Then} $\lb,\mu$ {\it are constant on } $M$.  {\it The distributions
} $D_{\lb},\ D_{\mu}$ {\it are both integrable if and only if}
$\n S=0$.
\par
\medskip
{\it Proof.} Let us note first that $p=\dim \ker(S-\lb I),
q=\dim \ker(S-\mu I)$ are constant on $M$ as $M_S=M$. We have also
$p\lb+q\mu=$tr$S$ and tr$S$ is constant on $M$. Hence
      $$p\n\lb+q\n\mu=0\tag 2.1$$
on $M$. Note that  $\la D_\lb,\n\lb\ra=0,\la D_\mu,\n\mu\ra=0$
thus $\la\n\lb,\n\mu\ra=0$ and from (2.1) it follows that $\n\lb=\n\mu=0$.
Let us note further that if $D_{\lb}$ is integrable then $\n S(X,Y)=0$
and $\n_{X}Y\in D_{\mu}$ if $X\in D_{\lb}$ and $Y\in D_{\mu}$.
We have $\n S(X,Y)=(\mu I-S)(\n_XY) \in D_{\lb}$ as $D_{\lb}$ is
orthogonal to $D_{\mu}$.  Let $Z\in \G(D_{\lb})$ then for any
$X\in \G(D_{\lb}),\ Y\in \frak X(M)$
we have
$$ \la\n S(X,Y),Z\ra=\la Y,\n S(X,Z)\ra=0$$
as $\n S(X,Z)=0$ ($D_{\lb}$ is integrable!). Hence $\n S(X,Y)=0$
and $\n_{X}Y\in D_{\mu}$ if $X\in D_{\lb}$ and $Y\in D_{\mu}$. If
$D_{\mu}$ is also integrable then in view of Corrolary 1.4 $\n S=0$.
$\k$
\par
\medskip

We have also proved:
\par
\medskip
{\it Corollary 2.2.} Let $S$ be an $\Cal A$-tensor on $(M,g)$ with two
constant eigenvalues $\lb,\mu$. If $D_{\lb}$ is integrable then
$\n S(X,Y)=0$ for all $X\in\G(D_{\lb}),\ Y\in\G(D_{\mu}).$

\par
\medskip
{\it Corollary 2.3.} Let $(M,g)$ be an $\A$-manifold, whose Ricci tensor $S$
has exactly two eigenvalues $\lb,\mu$. Then $\lb,\mu$ are constant.
\par
\medskip
{\it Proof.} It is well known that if $(M,g)$ is an $\A$-manifold then
$S$ has constant trace tr$S=\tau$ (see [7]).$\k$
\par
\medskip
From now on we shall investigate $\A$-tensors with two constant
eigenvalues $\lb,\mu$ satisfying an additional condition
$\dim D_{\lb}=1.$  It follows that $D_{\lb}$ is integrable.
We shall assume also that $\1$ is orientable (it happens for example if
$\pi_1(M)$ does not include subgroups of an index 2). In the other case,
we may consider a manifold $(\bar M,\bar g)$ and $\A$-tensor $\bar S$
on $\bar M$  such that there exists a two-fold Riemannian covering
$p:\bar M\rightarrow M$ for which $dp\circ \bar S=S\circ dp$
and $\bar \1=ker(\bar S-\lb I)$ is orientable.
Let $\xi\in\G(\1)$ be a global section of $\1$ such that $\la\xi,\xi\ra=1$.
Then we have:
\par
\medskip
{\it Lemma 2.4.} The section $\xi$ is a Killing vector field on $(M,g)$
and $L_{\xi}(\G(\2))\subset\G(\2)$ (i.e. for all $X\in \G(\2)\ L_{\xi}X=
[\xi,X]\in \G(\2)$.)
\par
\medskip
{\it Proof.} Let us denote by $T$ an  endomorphism of $TM$ defined by
$TX=\n_X\xi$. If $\Phi(X,Y)=\la SX,Y\ra$ then $\Phi(\xi,X)=\lb\la\xi,X\ra.$ Hence
$$ \n \Phi(Y,\xi,X)+\Phi(TY,X)=\lb\la TY,X\ra. \tag 2.2$$
Let us take $X=Y\in\2$ in (2.2). Since $\n S(X,X)=0$ ($\mu$ is constant)
we obtain $\Phi(TX,X)=\lb\la TX,X\ra.$
On the other hand $SX=\mu X$. Consequently, $\Phi(TX,X)=\mu \la TX,X\ra.$
Hence
$$ \la TX,X\ra=0 ,\  X \in \2 .  \tag 2.3$$
Since $\lb$ is constant, we have $\n_{\xi}\xi \in \1$. But we also have
$\la\n_{\xi}\xi,\xi\ra=0$ in view of $\p \xi\p=1.$  Hence $T\xi=\n_{\xi}\xi=0$.
 From (2.3) we obtain for all $X,Y\in \2$
 $$ \la TX,Y\ra=-\la X,TY\ra.  \tag 2.4$$
 From $\la\xi,\xi\ra=1$ it follows
 $$\la \n_X\xi,\xi\ra=\la TX,\xi\ra=0 \tag 2.5$$
 for all $X\in TM$, hence im$T\subset \2$.
 Consequently $\la T\xi,X \ra=-\la TX,\xi\ra=0$.
 Therefore  $\la TX,Y\ra=-\la X,TY\ra$ for all $X,Y\in TM$.
 Since $A_{\xi}=L_{\xi} -\n_{\xi}=-T$ is an antysymmetric operator
 it follows that $\xi$ is a Killing vector field, $L_{\xi}g=0$.
 If $X\in\2$ then $\la X,\xi\ra=0$ and consequently $0=L_{\xi}(\la X,\xi\ra)=
 \la [X,\xi],\xi\ra.$ $\k$
 \par
 \medskip
 Thus we obtain for the curvature tensor $R$ of $(M,g)$ (see [11]):
\par
\medskip
{\it Lemma 2.5.} Under the above assumptions
the following relations hold for all $X,Y\in TM$:
$$(a)\ \ \  R(X,\xi)Y=\n T(X,Y);\ (b)\ \ \n T(X,\xi)=-T^2X. \tag 2.6$$
 In particular
$$\rho (\xi,X)=-\la X,tr_g \n T \ra.\tag 2.7$$
\par
\medskip
{\it Proof.} The above formulas are valid for any Killing vector field
$\xi$, (2.6a) is well known. From $T\xi=0$ it follows
$\n T(X,\xi)+T(\n_X\xi)=0$ hence $\n T(X,\xi)=-T(TX)=-T^2X.$
From (a) we obtain
$$\rho(\xi,X)=\sum_{i=1}^n \la \n T(E_i,X),E_i\ra=-\la X,tr_g \n T\ra$$
where $\{E_1,E_2,..,E_n\}$ is a local orthonormal frame on $M.$ $\k$
\par
\medskip
{\it Corollary 2.5.} For all $X,Y\in TM$ we have:
$$ R(X,\xi)\xi=-T^2X,\ \la R(X,\xi)\xi,Y\ra=-\la T^2X,Y\ra=\la TX,TY\ra.$$
\par
\medskip
{\it Corollary 2.6.} Let $(M,g)$ be an $\A$-manifold whose Ricci
endomorphism $S$
has exactly two eigenvalues $\lb,\mu$ and such that dim$\1=1.$
Let $\xi$ be a local section of $\1$ and $\la\xi,\xi\ra=1.$
Then for all $X\in TM$:
$$ \la R(X,\xi)\xi,X\ra=\p TX\p^2 \tag 2.8$$
and also $\lb=\p T\p^2\geq 0$, tr$_g\n T=-\lb \xi.$

\par
\medskip

{\it Proof.} Let us note that if $\rho$ is a Ricci tensor of $(M,g)$
and $\{E_1,E_2,..,E_n\}$ is an orthonormal local frame on $M,\ \dim M=n$
 then $$\lb=\rho(\xi,\xi)=\sum_{k=1}^n \la R(E_k,\xi)\xi,E_k\ra=
\sum_{k=1}^n \la TE_k,TE_k\ra=\p T\p^2.$$
From Lemma 2.5 it follows $\rho(\xi,X)=-\la X,$tr$_g\n T\ra$ hence
tr$_g\n T=-\lb \xi.$
\par
\medskip
{\it Lemma 2.7.} Let $S,T$ be as above and define a 2-form
$\Omega\in\A ^2(M)$ by $\Omega(X,Y)=\la TX,Y\ra.$ Then $\Omega=d\th$ where
$\th (X)=\la \xi,X\ra.$ Consequently $\Omega$ is closed.
\par
\medskip
{\it Proof.} It follows from the simple computation:
$$d\th(X,Y)=\n_X\th(Y)-\n_Y\th(X)=
\la TX,Y\ra-\la TY,X\ra=2\la TX,Y\ra.$$ $\k$

Let us note that div$\xi=$tr$T=0$ hence $\th$ is co-closed.
\par
\medskip
{\bf Theorem 2.8.} {\it Let} $S$ {\it be an} $\A${\it -tensor
with two constant different eigenvalues } $\lb,\mu$ {\it and the distribution}
$\1$ {\it be one dimensional. Then the distribution } $\2$ {\it
is integrable if and only if } $\Omega=0$ {\it which means that }
$\th$ {\it is closed. The tensor} $S$ {\it is
parallel if and only if}  $\Omega=0$.
\par
\medskip
{\it Proof.} If $X,Y\in \G(\2)$ then
$\Omega(X,Y)=d\th(X,Y)=X\th(Y)-Y\th(X)-\th([X,Y]).$
Thus
$$\th([X,Y])=-\Omega(X,Y).$$
Hence $[X,Y]\in\G(\2)$ if and only if $\Omega=d\th=0.$ $\k$

\par
\medskip
{\it Corollary 2.9.} If $M$ is compact and $\lb\ne\mu$ then $S$
 is parallel if and
 only if the form $\th$ is harmonic and then  $b_1(M)>0.$

\par
\medskip
{\bf Theorem 2.10.} {\it Let } $(M,g)$ {\it be a Riemannian manifold
and let} $\xi$ {\it be a Killing vector field on}  $M$ {\it such that}
 $\la \xi,\xi\ra=1.$ {\it  Let us define a  tensor} $S$ {\it as follows:}
  $S\xi=\lb\xi$
 {\it and} $SX=\mu X \ \ if \la X,\xi\ra=0$ {\it where } $\lb,\mu$ {\it  are two
 different real numbers. Then} $S\in \A .$ {\it Moreover } $\n S=0$ {\it if
 and only if } $\n \xi=0.$
 \par
 \medskip
 {\it Proof.} We start with the lemma:
 \par
 \medskip
 {\it Lemma 2.11.} Let $S$ be a symmetric $(1,1)$-tensor ($\la SX,Y\ra
 =\la SY,X\ra$)
 with two constant eigenvalues $\lb,\mu$. Then $S\in \A$ if and only
 if the condition
 $$\n S(X,X)=0 \tag *$$
 is satisfied for all $X\in \1$ and all $X\in \2$.
\par
\medskip

 {\it Proof.} From (*) we obtain $\n S(X,Y)=-\n S(Y,X)$ if $X,Y\in \1$
 (resp. if $X,Y\in \2$). If $X\in \1$ (resp. if $X\in\2$) then
 $$ \n S(Y,X)=(\lb I-S)(\n_Y X) \in \2 (\ \text{resp.}
 (\mu I-S)(\n_Y X)\in \1). \tag 2.9$$

 Thus $\frak C_{X,Y,Z}\la\n S(X,Y),Z\ra=0$ if $X,Y,Z\in \1 (\2)$ where
 $\frak C$ denotes the cyclic sum.
 Hence it is enough to prove that
 $$ \la\n S(X,Y),Z\ra+\la\n S(Z,X),Y\ra+\la\n S(Y,Z),X\ra=0 \tag 2.10$$
 if $X,Y\in \1,Z\in\2$.
 From (2.9) it follows that
$ \la\n S(Z,X),Y\ra=0.$  We also have (in view of $\n  S(X,Y)=-\n S(Y,X)$)
 $$\la\n S(X,Y),Z\ra=-\la\n S(Y,X),Z\ra=-\la X,\n S(Y,Z)\ra.$$
 Hence (2.10) holds. Analogously we prove the case $X,Y\in \2,Z\in\1$.
 Thus we have proved that (2.9) holds if each $X,Y,Z$ belongs to one of
 the distributions $\1,\2$. Consequently, (2.10) holds for all
 $X,Y,Z\in TM.$$\k$
 \par
 \medskip

 Now we shall complete the proof of the theorem. Let us note that
 $\n S(\xi,\xi)=(\lb I-S)(\n_{\xi}\xi)=0$ as $\n_{\xi}\xi=0$ ($\xi$ is
 a unit Killing field see [11],[15]).
  We shall show that $\n S(X,X)=0$
 for all $X\in \2$. If $X\in \G(\2)$ then $\la X,\xi\ra=0$ and
 $\la\n_XX,\xi\ra+\la X,TX\ra=0$, where $TX=\n_X\xi.$ But $\la X,TX\ra=0$
 since $\xi$ is a Killing vector field. Thus, $\n_XX\in \G(\2).$
Consequently $\n S(X,X)=(\mu I-S)(\n_XX)=0.$ Hence our theorem
follows from Lemma 2.11.$\k$
\par
\medskip

{\it Remark 2.12.} Let us note that if a space $(M,g')$ admits a
nonvanishing Killing vector field $\xi$  then there exists
a metric $g$ on $M$ conformally equivalent to $g'$
such that $\xi$ is a unit Killing vector field of $(M,g)$
(see [15]). Hence on every such manifold there exists a Riemannian
structure admiting an $\A$-tensor
with two different eigenvalues.
In particular if $M$ admits an effective free action
of the group $S^1$ then it admits a Riemannian metric
$g$  and an $\A$-tensor $S$ on $(M,g)$ with two different
eigenvalues such that  the fundamental vector field $\xi$ of the
action of the group $S^1$ is an eigenfield of $S$ (see [W]).

\par
\bigskip
{\bf   3.  The structure of $\A$-manifold on a $S^1$-principal
fibre bundle.} In this section we shall construct the Riemannian
metric $g$ on the $S^1$-principal fibre bundle $P$ over
a K\"ahler-Einstein manifold $(M,g_*,J)$ such that $(P,g)$ is an
$\A$-manifold.  We generalize  A.Gray's construction of
the $\A$-structure on the $S^1$-bundle $P=S^3$ over $M=S^2.$

Let $(M,g_*)$ be a Riemannian space, $\dim M=m$, and let
$p:P\rightarrow M$ be a principal
$S^1$-bundle over $M$. Let $\bt\in \A^1(P)$ be a connection form on $P$
and let $\bar{\xi}$ be the fundamental vector field of the action of
the group $S^1$ on $P$. Thus, $\bt(\bar\xi)=1$. By $\bO$ we denote
the curvature form $d\bt$. Let us define for a number $c>0$
the metric $g=g_c$ on $P$ as follows:
$$g_c(X,Y)=c^2\bt (X)\bt (Y)+g_*(dp(X),dp(Y)). \tag g$$
Then $p:(P,g_c)\rightarrow (M,g)$ is a Riemannian submersion. Let us note
that $L_{\xi}\th=0$ where we define $\th=c\bt$
and $\xi=\frac1c\bar{\xi}$ ($\bt$ is a connection form!). Hence
$L_{\xi}g_c=0$ which means that $\xi$ is a unit Killing vector field on
$(P,g_c)$. The form $\Omega=c\bO$ is projectable and let $p^*\wO=\Omega$,
where $\wO\in \A^2(M)$. Notice that $\th(X)=\la \xi,X\ra$ and
$\0=d\th=2\la TX,Y\ra$
where $T=\n \xi.$ The tensor $T$ satisfies relations:
$T\xi=0$ ($\n_{\xi}\xi=0$ since $\xi$ is a unit Killing vector field ) and
$L_{\xi}T=0$ . In fact
$$0=d(L_{\xi}\th)(X,Y))=L_{\xi}(\0)(X,Y)=2\la(L_{\xi}T)X,Y\ra$$
hence $L_{\xi}T=0$. It follows that there exists a tensor $\wT$ on $M$
such that $\wT\circ dp=dp\circ T$ and $\wO (X,Y)= 2\la\wT X,Y\ra_*$.
We shall check under what conditions $\xi$ is an eigenfield of the Ricci
tensor $S$ of $(P,g_c)$. Let us note that if $\wn$ is a Levi-Civita
connection for $(M,g_*)$ then $(\wn_XY)^*=\n_{X^*}Y^*
-\frac12 \Cal V [X^*,Y^*]$
and $\la\wn_XY,Z\ra=\la\n_{X^*}Y^*,Z^*\ra$
where $X,Y,Z\in \frak X(M)$, $*$ denotes the horizontal lift and
$\Cal VX$ (respectively $\Cal HX$) denotes a vertical (horizontal)
part of $X$.
Hence  if $\{E_1,E_2,..,E_{m}\}$ is an orthonormal local frame on $M$
then we have:
$$\gather
\la tr_{g_*}\wn\wT,Y\ra_*=\sum_{i=1}^{m} \la\wn \wT (E_i,E_i),Y\ra_*
=\sum_{i=1}^{m}(\la\wn_{E_i}(\wT E_i),Y\ra_*\\-\la\wT(\wn_{E_i}E_i,Y\ra_*)=
\sum_{i=1}^{m}(\la\n_{E_i^*}(T E_i^*),Y^*\ra-\la T(\n_{E_i^*}E_i^*,Y^*\ra)\\=
\la tr_g \n T,Y^*\ra.
\endgather$$
Consequently,
$$ \delta \wO=-\sum_{i=1}^{m}\wn_{E_i}\wO(E_i,Y)=-2\la tr\wn\wT,Y\ra_*=
-2\la tr\n T,Y^*\ra\tag A$$
as $\n T(\xi,\xi)=0$. From (2.7) it follows that $\xi$ is an eigenfield
of the Ricci tensor if and only if $\delta \wO=0$ which is equivalent to
$tr_g\n T\p \xi$. In view of a relation $p^*\wO=d\th$ the form $\wO$ is
closed. Hence we obtain:
\par
\medskip
{\bf Proposition 3.1.} {\it Let} $p:P\rightarrow M$ {\it be a }
$S^1${\it-principal fibre bundle over} $(M,g_*)$ {\it and}
 $\bt$ {\it be a connection form on P. If we define the metric g on
 P by the formula (g) then the fundamental vector field} $\bx=c\xi$
 {\it  of the action of the group } $S^1$ {\it is a Killing vector field.
 The field } $\xi$ {\it is the eigenfield of the Ricci tensor S of} $(P,g)$
 {\it i.e. }$S\xi=\lb\xi$ {\it if and only if} $\delta \wO=0$
 {\it and then} $\lb=\p\wT\p^2,\ S\xi=\p\wT\p^2\xi.$
 {\it If M is compact then the above condition means that } $\wO$ {\it
 is harmonic} $(\Delta\wO=0).$
\par
\medskip
{\it Proof.} The proposition follows from (2.7) and (2.8) and the relation
$\p\wT\p=\p T\p.$ $\k$
\par
\medskip
{\it Remark 3.2.} Let us note that $\lb=\rho(\xi,\xi)=\p\n\xi\p^2=
\p T\p^2=\frac14\p\0\p^2$ hence $\lb$ is constant if and only if $\0$ has
constant length.
\par
\medskip
Now we find under what conditions the Ricci tensor $S$ satisfies the
relation $SX^*=\mu X^*$ for all $X\in\frak X(M)$. We shall use
O'Neill formulas (see [13],[1]). The fibers of $p:P\rightarrow M$ are
totally geodesic ($\n_{\xi}\xi=0$!.)  Hence the O'Neills tensor $T$
vanishes. We shall compute the tensor $A$:
$$A_EF=\Cal V(\n_{\Cal H E}\Cal HF)+ \Cal H(\n_{\Cal H E}\Cal VF).$$
It is easy to check that:
$$A_EF=\la E,TF\ra\xi+\la\xi,F\ra TE. \tag 3.1$$
Consequently, if $U,V$ are horizontal vectors then
$$\p A_UV\p^2=\la V,TU\ra^2\tag 3.2$$
and $A_UV=-\frac12\0(U,V)\xi.$ Hence
$$K(U\w V)=K_*(U_*\w V_*)-3\la U,TV\ra^2 \tag 3.3$$
and
$$K(U\w\xi)=\p TU\p^2. \tag 3.4$$
where $K$ denotes a sectional curvature and $ (U_*)^*=U,\ (V_*)^*=U$.
Let $\{ E_0=\xi,E_1,E_2,..,E_{m}\}$ be a local orthonormal frame
on $P$. Then we obtain for the Ricci tensor $\rho$ the formula ($U$ is a
unit horizontal vector):
$$\rho(U,U)=K(U\w\xi)+\sum_{i=1}^{m}K(U\wedge E_i)=\rho_*(U_*,U_*)
-2\p TU\p^2.\tag 3.5$$
We also obtain a formula for the scalar curvature $\tau$ of $(P,g)$:
$$\tau=\tau_*-2\p T\p^2+\p T\p^2=\tau_*-\p T\p^2$$
where $\rho_*,\tau_*$ are the Ricci tensor and the scalar curvature of
$(M,g_*)$. Hence $SX^*=\mu X^*$ if and only if
$\rho_*(U_*,U_*)-2\p \wT U\p_*^2 = \mu\p U_*\p_*^2$ or equivalently if
$S_*+2\wT^2=\mu Id,$ where $\rho_*(X,Y)=\la S_*X,Y\ra_*.$
If $(M,g_*)$ is a Riemannian space then an integral 2-form
$\0$ (i.e. $\{\0\}\in H^2(M,\Bbb{Z})$) determines a $S^1$
principal fibre bundle $p:P\rightarrow M$ and a connection form
$\bt\in\A^1(P)$ such that $d\bt=2\pi p^*\0.$  We can construct on $P$
a Riemannian metric $g_c$ using formula (g) and then $p$ is a Riemannian
submersion. The fundamental field $\xi$ is an eigenfield of the Ricci
tensor of $(P,g_c)$ if and only if $\delta \0=0.$
We shall prove a theorem similar to the results in [1] (we use the
above notation):
\par
\medskip
{\bf Theorem 3.3.} {\it Let} $(M,g_*,J)$ {\it be a
K\"ahler-Einstein manifold with nonvanishing scalar curvature} $\tau_*$.
{\it Then there exists a } $S^1${\it-principal fibre bundle }
$p:P\rightarrow M$ {\it and a connection form } $\bt$ {\it on P such that
} $d\bt=-\a p^*\omega$ {\it where } $\omega$ {\it is
a K\"ahler form of } $(M,g_*,J)$ {\it and} $\a=\frac{\tau_*}{2n},\
 n=\dim_{\Bbb{C}}M.$ {\it If for } $c>0$ {\it we define}
 $\ g_c=c^2\bt\otimes\bt+p^*g_*$
 {\it then} $\wT=- \frac12 c\a J$ {\it and } $\delta\wO=0$.
 {\it Consequently,}
 $(P,g_c)$ {\it is an } $\A${\it-manifold and for } $c^2\ne\frac2{(n+1)\a}$
 {\it is a proper} $\A${\it-manifold with two constant eigenvalues}
 $\lb=\frac12nc^2\a^2$ {\it and} $\mu=\a(1-\frac12\a c^2).$
\par
\medskip
{\it Proof.} (See also Th.9.76 in [1].) Let $P$ be a $S^1$-bundle
determined (see [8]) by the  first Chern class $c_1(M)$ of $(M,g_*,J)$.
Let us recall that $c_1(M)=\{-\frac{\rho_J}{2\pi}\}$, where
$\rho_J$ is the Ricci form of $(M,g_*,J)$ i.e. $\rho_J(X,Y)=\rho_*(X,JY).$
As $(M,g_*)$ is an Einstein space we have
$\rho_*(X,Y)=\frac{\tau_*}{2n} \la X,Y\ra_*$ where $\tau_*$ denotes the scalar
curvature of $(M,g_*)$. Let $\bt$ be a connection form on $P$ such that
$d\bt=-\a p^*\omega$ where $\omega$ is the K\"ahler form of $(M,g_*,J)$
($\omega(X,Y)=\la X,JY\ra$) and $\a=\frac{\tau_*}{2n}$. Note that $P$ can be
realized as a subbundle of the anti-canonical line bundle
$K^*=\bigwedge^nT^{(1,0)}M$. If we define $P=\{\g\in K^*:\la\g,\bar\g\ra=1\}$
then $P$ with the induced metric connection satisfies the above conditions.
Let us define a metric $g$ on $P$ by the formula (g) then a form
$\0=cd\bt$ satisfies the relation $\0=-c\a p^*\omega$
thus $\wO=-c\a\omega$
($\0,\wO,\wT,T$ are defined as in Prop.3.1).
Hence $2\wT=-c\a J$ and $\p T\p^2= \frac12nc^2\a^2.$ It is clear that
$\delta\wO=0$. Consequently, $\xi$ is an eigenfield of the Ricci tensor $S$
of $(P,g)$ and
$$S\xi=\frac12 nc^2\a^2\xi. \tag 3.6$$
From (3.5) it follows that
$$\rho(U,U)=\a\la U,U\ra-\frac12 c^2\a^2\la U,U\ra
=\a(1-\frac12 c^2\a)\la U,U\ra \tag 3.7$$
for all horizontal vectors $U\in TP$. From (3.6) it follows that $S$
preserves the distribution $H$ of horizontal vectors of $P$ i.e.
$SH\subset H$, hence from (3.7) it is clear that
$$SU=\a(1-\frac12 c^2\a)U\tag 3.8$$
for any $U\in H$. Let us denote $\lb=\frac12 nc^2\a^2,\
\mu=\a(1-\frac12 c^2\a)$.
If $c^2=\frac2{(n+1)\a}$ then $\lb=\mu$ and $(P,g_c)$ is an Einstein space.
Let us assume that $c^2\ne\frac2{(n+1)\a}$. Then the assumptions of
 Th.2.10 are satisfied hence $S\in\A$ which means that $(P,g_c)$ is an
 $\A$-manifold. The space $(P,g_c)$ is a proper $\A$-space if
 $c^2\ne\frac2{(n+1)\a}$ as we have assumed that $\tau_*\ne0$
 hence $\wO\ne 0$ (see Th.2.8)$\k$
\par
\bigskip

{\it Remark 3.4.} Notice that if $(M,g_*,J)$ is a compact K\"ahler-Einstein
manifold then $(P,g_c)$ is a compact $\A$-manifold.
\par
\medskip
{\it Corollary 3.5.} If $(M,g_*)$ is a closed Riemannian surface
of constant non-zero curvature $K\in \Bbb{R}$ then there exists
a $S^1$-bundle $P$ over $M$ and a family of Riemannian structures $(P,g_c)$
 ($c>0,c^2\ne\frac1{K}$) on $P$
 such that $(P,g_c)$ is a  proper compact $\A$-manifold and a submersion
 $p:P\rightarrow M$ is a Riemannian submersion. If $K>0$ then we obtain
 A. Gray's examples (see [7] p.267).
\par
\medskip
Let us recall that N.Koiso and Y.Sakane have constructed explicit
 locally non-homogene\-ous
examples of compact K\"ahler-Einstein manifolds with arbitrary cohomogeneity
 $d\in \Bbb{N}$ (see [9], [10]) and E.Calabi has constructed
 locally non-homogene\-ous complete K\"ahler-Einstein manifolds for
 every dimension $2n,\ \ n>1$.
 Hence we have:
 \par
 \bigskip
{\it Corollary 3.6.} If $(M,g_*,J)$ is a compact non-homoge\-nous
K\"ahler-Einstein
manifold with $\tau_*\ne 0$ and cohomogeneity $d$ (see  [9],[10])
then the space $(P,g_c)$ (for $c$ satisfying a condition
$c^2\ne \frac{4n}{(n+1)\tau_*}$)
is a non-homogon\-eous
proper compact $\A$-manifold of cohomogeneity  $d$.
\par
\medskip

{\it Proof.} Let us note that if $X$ is a Killing vector field on
$(P,g_c)$ then $L_XS=0$ and  from (3.6) we obtain
$S([X,\xi])=\lb [X,\xi]$
thus,  (we assume that $\lb\ne\mu$), $[X,\xi]\parallel \xi.$ On the other
hand, the relation
$<\xi,\xi>=1$ yields $<[X,\xi],\xi>=0$. Hence $[X,\xi]=0$. It follows that
$X$ is projectable $\Cal HX=X_1^*$ where $X_1\in \frak X(M)$ is a Killing
vector field on $M$. Consequently every Killing vector field on $P$
is projectable  and cohomg$(P)=$cohomg$(M)$.$\k$
\par
\medskip
{\it Corrolary 3.7.} If $(M,g_*,J)\ \ $ is a locally non-homogeneous
complete K\"ahler-Einstein
manifold with negative scalar curvature and with dim$M=2n>2$ (see [2],[16])
  then $(P,g_c)$
 is a complete locally non-homogeneous
$\A$-manifold of dimension $2n+1$ giving an answer to the
open problem in [1] (p.451,16.56 (i)). Hence for every odd number
$m>3$ we have constructed locally non-homogeneous proper $\A$-manifold
$(M,g)$ with dim$M=m$.  Let us note that there are many such
 compact manifolds given by Calabi-Yau theorem, hence we obtain many
 examples of non-homogeneous proper compact $\A$-manifolds (however they
 are not given explicitly).
 \par
 \medskip
 {\it Remark 3.8} Let us note that for $c=\frac2{\a}=\frac{4n}{\tau_*}$ the
 manifold $P$ is a Sasakian $\A$-manifold. We classify Sasakian
 $\A$-manifolds in the forthcoming paper.

\par
\bigskip

\centerline{\bf References.}
\par
\medskip
\cite{1} A.L. Besse, `Einstein Manifolds'
{\it Springer Verlag} Berlin Heidelberg (1987)
\par
\medskip
\cite{2} E.Calabi, {\it A construction of nonhomogeneous Einstein
metrics} Proc. Sympos. Pure Math., vol {\bf 27}, Part II, Amer. Math. Soc.,
Providence, R.I., 1975, 17-24.

\par
\medskip
\cite{3} J.E. D'Atri, {\it Geodesic spheres and symmetries in naturally
 reductive homogeneous spaces} Mich. Math. J. {\bf 22} (1975), 71-76.
 \par
 \medskip
 \cite{4} J.E. D'Atri and H.K. Nickerson, {\it Divergence preserving
 geodesic symmetries} J. Diff.Geom {\bf 3} (1969), 467-476.
 \par
 \medskip
 \cite{5} J.E. D'Atri and H.K. Nickerson, {\it Geodesic symmetries
 in spaces with special curvature tensor} J. Diff.Geom {\bf 9} (1974),
 251-262.

\par
\medskip
\cite{6}  A. Derdzi\'nski, {\it Classification of certain compact
Riemannian manifolds with harmonic curvature and non-parallel
Ricci tensor} Math. Z. {\bf 172}  (1980), 273-280.
\par
\medskip
\cite{7} A. Gray, {\it Einstein-like manifolds which are not Einstein}
Geom. Dedicata {\bf 7} (1978), 259-280.
\par
\medskip
\cite{8}  S. Kobayashi, {\it Principal fibre bundles with the
1-dimensional toroidal group} T\^ohoku Math. J.{\bf 8} (1956), 29-45.
\par
\medskip
\cite{9} N. Koiso and Y. Sakane, {\it Non-homogeneous K\"ahler-Einstein
metrics on compact complex manifolds}, Curvature and Topology of
Riemannian manifolds, Proceedings 1985, Lecture Notes in Math.1201,
Springer-Verlag, (1986),165-179.
\par
\medskip
\cite{10} N. Koiso and Y. Sakane,{\it Non-homogeneous K\"ahler-Einstein
metrics on compact complex manifolds II}, Osaka J. Math. {\bf 25}
 (1988),  933-959.
\par
\medskip
\cite{11} S.Kobayashi and K.Nomizu,  Foundations of differential
geometry. vol I, {\it Interscience Publishers}, New York, London
(1963).
\par
\medskip
\cite{12} O. Kowalski and L. Vanhecke, {\it Riemannian Manifolds
with Homogeneous Geodesics} Bolletino U.M.I. 7, {\bf 5}-B (1991),
189-246.
\par
\medskip
\cite{13} B. O'Neill,  {\it The fundamental equations of a submersion}
Mich. Math. J. {\bf 13} (1966), 459-469.

\par
\medskip
\cite{14} K. Sekigawa and L. Vanhecke, {\it Symplectic geodesic symmetries
on K\"ahler manifolds} Quart.J.Math Oxford (2), {\bf 37} (1986), 95-103.
\par
\medskip
\cite{15} A.W. Wadsley, {\it Geodesic foliations by circles}
J. Diff. Geom. {\bf 10}  (1975), 541 - 549.
\par
\medskip
\cite{16} J.A.Wolf, {\it On Calabi's Inhomogeneous Einstein-K\"ahler
manifolds } Proc. Amer. Math. Soc, {\bf 63}, No 2, 1977, 287-288.

\par
\medskip
\noindent Institute of Mathematics
\par
\noindent Technical University of Cracow
\par
\noindent Warszawska 24
\par
\noindent 31-155 Krak/ow, POLAND.

\enddocument